\documentclass{article}

\usepackage{amsmath,amssymb,amsthm}

\usepackage[T1]{fontenc}
\usepackage{libertine}
\usepackage[libertine]{newtxmath}

\numberwithin{equation}{section}
\newtheorem{thm}{Theorem}[section]
\newtheorem{cor}[thm]{Corollary}
\newtheorem{prop}[thm]{Proposition}
\newtheorem{lem}[thm]{Lemma}

\theoremstyle{definition}
\newtheorem{defn}[thm]{Definition}
\theoremstyle{remark}
\newtheorem{rmk}[thm]{Remark}
\newtheorem{exam}[thm]{Example}

\usepackage{microtype}

\usepackage[all]{xy}

\usepackage{todonotes}

\usepackage{hyperref}
\hypersetup{
  colorlinks   = true, %Colours links instead of ugly boxes
  urlcolor     = blue, %Colour for external hyperlinks
  linkcolor    = blue, %Colour of internal links
  citecolor   = red %Colour of citations
}

\newcommand{\co}{\colon\thinspace}
\newcommand{\mb}[1]{\mathbb{#1}}

\newcommand{\op}{{op}}
\newcommand{\act}{{act}}
\newcommand{\too}{\xrightarrow}

\newcommand{\cC}{\mathcal{C}}
\newcommand{\dD}{\mathcal{D}}
\newcommand{\oO}{\mathcal{O}}

\newcommand{\wt}[1]{\widetilde{#1}}

\newcommand{\tH}{\widehat{H}}

\DeclareMathOperator{\Hom}{Hom}
\DeclareMathOperator{\Map}{Map}

\DeclareMathOperator{\Fun}{Fun}

\DeclareMathOperator{\Sp}{Sp}
\DeclareMathOperator{\Fin}{Fin}
\DeclareMathOperator{\Spaces}{\mathcal{S}}
\DeclareMathOperator{\Ch}{Ch}
\DeclareMathOperator{\Tate}{Tate}
\DeclareMathOperator{\Env}{Env}
\DeclareMathOperator{\Alg}{Alg}
\DeclareMathOperator{\LMod}{LMod}

\DeclareMathOperator{\Free}{Free}

\DeclareMathOperator{\Cat}{Cat}

\DeclareMathOperator*{\hocolim}{hocolim}

\DeclareMathOperator{\THH}{THH}

\DeclareMathOperator{\HH}{HH}
\DeclareMathOperator{\Tr}{Tr}
\DeclareMathOperator{\TBar}{Bar}

\title{Unwinding the relative Tate diagonal}
\author{Tyler Lawson\thanks{The author was partially supported by NSF
    grant 1560699.}}

\begin{document}
\maketitle

\begin{abstract}
  We show that a spectral sequence developed by Lipshitz and Treumann,
  for application to Heegaard Floer theory, converges to a localized
  form of topological Hochschild homology with coefficients. This
  allows us to show that the target of this spectral sequence can be
  identified with Hochschild homology when the topological
  Hochschild homology is torsion-free as a module over
  $\THH_*(\mb F_2)$, parallel to results of Mathew on degeneration of
  the Hodge-to-de Rham spectral sequence.

  To carry this out, we apply work of Nikolaus--Scholze to develop
  a general Tate diagonal for Hochschild-like diagrams of spectra that
  respect a decomposition into tensor products. This allows us to
  discuss the extent to which there can be a Tate diagonal for relative
  topological Hochschild homology.
\end{abstract}

\section{Introduction}
\label{sec:introduction}

\subsection*{Hochschild homology and cyclic covers}

The primary goal of this paper is to understand a particular formality
condition on Hochschild homology. Motivated by Heegaard Floer
homology of double covers, Lipshitz and Treumann developed a
noncommutative version of the Hodge-to-de Rham spectral sequence with
coefficients in \cite{lipshitz-treumann-tate}. Given a homologically
smooth and proper differential graded algebra over $\mb F_2$ and a
bounded differential graded $A$-bimodule $M$, they give a spectral
sequence with $E_1$-term
\[
  E^1_* = HH_*(A; M \otimes_A^L M).
\]
The $E_2$-term is Tate cohomology for an action of the cyclic group
$C_2$ \cite[Theorem 4]{lipshitz-treumann-tate}---this spectral
sequence arises from the Tate construction for $C_2$ acting on the
Hochschild complex $HH(A; M \otimes_A^L M)$. When $M = A$, the
$d_1$-differential is trivial and the $d_2$-differential is the
$b$-operator of Connes, and this spectral sequence is related to the
Hodge to de Rham spectral sequence \cite{kaledin-hodgederham}. Under
an assumption called \emph{$\pi$-formality}, they show that their
spectral sequence converges to $HH_*(A;M)$. Without $\pi$-formality,
the target of their spectral sequence is not easy to identify, and
some of the steps in their identification make use of non-additive
maps $x \mapsto x \otimes x$ on the homology level that do not lift to
the chain level.

The main goal of this paper is to identify the target of
Lipshitz--Treumann's spectral sequence with a periodic variant of
B\"okstedt's topological Hochschild homology with coefficients
\cite{bokstedt-thh}. To begin stating our results, we recall the
following result of Bhatt--Morrow--Scholze
{\cite{bhatt-morrow-scholze-thh}} on topological Hochschild homology
of perfect rings---for the ground field $\mb F_p$ this is due to
B\"okstedt, and for perfect fields, such as finite fields, this result
was previously known by work of Hesselholt and Madsen
\cite{hesselholt-madsen-finitealgebras}. If $R$ is a perfect ring of
characteristic $p$, the topological Hochschild homology $\THH_*(R)$ is
a polynomial algebra $R[u]$ on a generator in degree 2. Moreover, the
Tate cohomology ring $\tH^{-*}(C_p;R)$ is an algebra over $\THH_*(R)$
whose underlying module is
\[
  R[u^{\pm 1}] \cdot \{1, v\}.
\]
Here $1$, $v$, and $u^{-1}$ are generators of $H^i(C_p;R)$ for
$i = 0$, $1$, and $2$ respectively.

\begin{thm}
  \label{thm:mainhomological}
  Suppose that $\mb F$ is a perfect field of characteristic $p$, $A$ is a
  homologically smooth differential graded $\mb F$-algebra, and $M$ is a
  bounded $A$-bimodule which is finitely generated over $\mb F$. Then
  there is a Tate cohomology spectral sequence with $E_2$-term
  \[
    \tH^*(C_p; HH^{\mb F}_*(A, \underbrace{M \otimes^L_A \dots
    \otimes^L_A M}_{p})) \Rightarrow \tH^*(C_p; \mb F) \otimes_{\mb
    F[u]}\THH_*(A;M).
  \]
\end{thm}

Both sides of this periodic: degree $d$ is canonically isomorphic to
degree $(d+2)$ for all $d$ (and $d+1$ if $p=2$). By combining this
with base-change results for topological Hochschild homology, we will
arrive at the following result.

\begin{thm}
  \label{thm:spectralsequence}
  Suppose, under the hypotheses of Theorem~\ref{thm:mainhomological},
  that $\THH_*(A;M)$ is torsion-free as a module over
  $\THH_*(\mb F) = \mb F[u]$. Then the spectral sequence of
  Theorem~\ref{thm:mainhomological} reduces to an ungraded spectral
  sequence of the form
  \[
    \tH^0(C_2; HH^{\mb F}_*(A, M \otimes^L_A M)) \Rightarrow HH^{\mb
      F}_*(A;M)
  \]
  if $p=2$, and a similar $\mb Z/2$-graded spectral sequence
  when $p$ is odd.
\end{thm}

In some sense, $\THH(A;M)$ contains a derived interpolation between
the values $u=0$ (where we get Hochschild homology and this Tate
spectral sequence) and $u=1$ (where we eliminate the grading), and
torsion-freeness allows us to identify the two
results. We suspect Lipshitz--Treumann's development of
$\pi$-formality gives conditions under which topological Hochschild
homology is torsion-free.

\subsection*{Algebra and homotopy theory}

The connection between algebra and stable homotopy theory, and in
particular topological Hochschild homology, arises through the
following translation procedure.
\begin{itemize}
\item For a commutative ring $R$, there is an Eilenberg--Mac Lane
  spectrum $HR$ which has the structure of a commutative algebra.
\item The category of differential graded $R$-modules is equivalent
  (in a derived sense) to the category of $HR$-modules: this
  translation takes a complex $V$ to a spectrum $HV$ such that
  $H_*(V) \cong \pi_* HV$.\footnote{This is a high-powered version of
    the \emph{Dold--Kan equivalence} between chain complexes and
    simplicial abelian groups.}
\item This equivalence is symmetric monoidal, in the sense of
  \cite{lurie-higheralgebra}, and it takes the (derived) tensor over
  $R$ to the tensor over $HR$, known as the relative smash product.
\item This equivalence preserves homotopy limits, homotopy colimits,
  and Tate constructions.
\end{itemize}

In particular, under this correspondence a differential graded
$R$-algebra $A$ lifts to an $HR$-algebra $HA$, an $R$-linear
differential graded $A$-bimodule lifts to an $HR$-linear
$HA$-bimodule, and the relative Hochschild complex $\HH^R(A;M)$ lifts
to topological Hochschild homology $\THH^{HR}(HA;HM)$. Moreover, there
is a base-change formula
\[
  \THH^{HR}(HA;HM) \simeq HR \otimes_{\THH(HR)} \THH(HA;HM)
\]
due to McCarthy and Minasian. Together these equivalences relate
ordinary Hochschild homology with a base-change of $\THH$ \cite[\S
5]{mccarthy-minasian-hkr}. A sketched discussion of the translation
procedure will occupy \S\ref{sec:translation}.

The connection to equivariant stable theory was observed by Kaledin
\cite{kaledin-hodgederham}, and our methods are very similar to those
of Mathew \cite{mathew-degeneration}: however, where Mathew makes use
of the circle action on $\THH(A)$, we make use of actions of cyclic
groups on $\THH$ with certain coefficients. In these terms,
Theorem~\ref{thm:mainhomological} is a consequence of the following.

\begin{thm}
  \label{thm:main}
  Suppose that $k$ is a commutative ring spectrum, $A$ is a
  $k$-algebra, and $M$ is a $k$-linear $A$-bimodule. Then there exists
  an action of $C_p$ on the topological Hochschild homology
  $\THH(A;M \otimes_A \dots \otimes_A M)$ and a natural
  \emph{relative $\THH$-diagonal}
  \[
    k^{tC_p} \otimes_{\THH(k)} \THH(A,M)
    \to \left[\THH^k(A, M \otimes_A M \otimes_A \dots \otimes_A
      M)\right]^{tC_p}.
  \]
  If $A$ is a smooth $k$-algebra and the underlying $k$-module of $M$
  is perfect, this map is an equivalence.
\end{thm}

Once the relative $\THH$-diagonal is set up, the proof in
\S\ref{sec:smooth-algebras} that this is an equivalence will be a
relatively formal thick subcategory argument.

%\todo{something about lifting algebras to the sphere}

\subsection*{The relative Tate diagonal}

The results of this paper rely on a piece of nonalgebraic structure:
the \emph{Tate diagonal}, which plays a prominent role in equivariant
stable homotopy theory. For a spectrum $X$, the Tate diagonal is a
natural map
\[
  \Delta\co X \to (X^{\otimes p})^{tC_p}
\]
that enjoys a great deal of structure. The Tate diagonal is lax
symmetric monoidal, it is natural, and it is impervious to the action
of the cyclic group $C_p$ on $X^{\otimes p}$. These properties are
concisely encoded by Nikolaus--Scholze's expression of functoriality
on a category of finite free $C_p$-sets
\cite[III.3.8]{scholze-nikolaus-tc}.

The lax symmetric monoidality of the Tate diagonal allows us to
construct a relative version, using the tensor product over a
commutative ring spectrum $k$ rather than over the sphere. For
individual modules, the relative Tate diagonal behaves very similarly
to the ordinary one. However, \emph{the functoriality of the relative
  diagonal is less strong}: the lax symmetric monoidal compatibility
does not play well with the cyclic group invariance. We would
like to extend Nikolaus--Scholze's functoriality for the Tate diagonal
in a way that exhibits the extent to which there can be a relative
version.

The starting point is the observation that the cyclic bar construction
$Z(A)$ that builds the Hochschild complex is not just a simplicial
spectrum $\Delta^\op \to \Sp$: there is a decomposition of each
simplicial degree into a formal smash product of factors, and the
structure maps respect this decomposition.  This lifts it to a
simplicial object in the \emph{symmetric monoidal envelope} $\Env(\Sp)$,
which we will discuss in \S\ref{sec:envelopes}.

Associated to a diagram $X\co K \to \Env(\cC)$ in a symmetric monoidal
envelope, there is an underlying diagram of finite sets $K \to Fin$
representing the decomposition into tensor factors: we will define the
\emph{shape} $|X|$ to be the resulting simplicial set. For example,
the shape of $\THH(A)$ or $\THH(A;M)$ will be the circle $S^1$. Given
a diagram $X\co K \to \Env(\cC)$ in a symmetric monoidal envelope with
shape $|X|$ and a principal $C_p$-bundle classified by a map
$f\co |X| \to BC_p$, we will associate a new diagram \emph{unwinding}
$X$: a $C_p$-equivariant diagram $\psi^f X\co K \to \Env(\cC)$. For
any $i \in K$, $\psi^f(i)$ is isomorphic to the iterated formal tensor
$X(i)^{\otimes p}$, but the maps in $K$ go to maps between tensor
powers that make use of the structure of the principal bundle.

Let $k^{\otimes Z}$ denote the \emph{Loday construction}
\cite{schlichtkrull-loday} (sometimes called the \emph{factorization
  homology} \cite{francis-tangentcomplex}) of $Z$ with coefficients in
$k$, often also denoted by $\mathcal{L}_Z(k)$, $\int_{Z} k$, or
$Z \otimes k$.
\begin{prop}
  Let $k$ be a commutative ring spectrum and $K$ a sifted index
  category. For a diagram $X\co K \to \Env(\LMod_k)$ in the symmetric
  monoidal envelope and a principal $C_p$-bundle $f\co |X| \to BC_p$
  over the shape of $X$, there is a natural map
  \[
    k^{tC_p} 
    \mathop\otimes_{k^{\otimes |X|}} \left(\hocolim_{i \in K} \bigotimes^{\Sp} X(i)\right) \to
    \left(\hocolim_{i \in K} \bigotimes^{\LMod_k} \psi^f X(i)\right)^{tC_p}
  \]
  called the \emph{relative Tate diagonal}.
\end{prop}

There is asymmetry between the tensor products in the source and
target of the relative Tate diagonal---the source tensor takes place
in spectra and the target tensor takes place in $k$-modules. We have
traced several of our own misunderstandings, including a mistaken
assertion that there are cyclotomic structures on relative $\THH$ and
monoidality properties of a Tate diagonal on relative $\THH$ with
coefficients, back to this root. The relative Tate diagonal \emph{does
  not} imply that there is a $k$-module Tate diagonal
\[
  \hocolim_{i \in K} \bigotimes^{\LMod_k} X(i) \to \left(\hocolim_{i \in K}
  \bigotimes^{\LMod_k} \psi^f X(i))\right)^{tC_p}
\]
unless the map $k^{\otimes |X|} \to k^{tC_p}$ factors through the
augmentation $k^{\otimes |X|} \to k$. This only holds in a few
circumstances, such as when we can fix a trivialization of the bundle
classified by $f$.\footnote{Thomas Nikolaus has pointed out to us that
  such a factorization through the augmentation is also possible when
  $k$ is the spherical group algebra of a \emph{discrete} abelian
  group. See Remark~\ref{rmk:krause-nikolaus}.}  This is true, for
example, when the index category $K$ is a singleton, which allows one
to construct a natural Tate diagonal
$M \to (M \otimes_k \dots \otimes_k M)^{tC_p}$ and a $k$-module
version of the Hill--Hopkins--Ravenel norm \cite{hhr-kervaire}.

\subsection*{Acknowledgements}

The author would like to thank
Clark Barwick,
Andrew Blumberg,
Teena Gerhardt,
Lars Hesselholt,
Michael Hill,
Robert Lipshitz,
Michael Mandell,
Denis Nardin,
Thomas Nikolaus,
and
David Treumann
for their assistance and forbearance through this paper's long period
of development. The author would also like to thank the Max Planck
Institute for Mathematics in Bonn for their hospitality and financial
support while this paper was written.

\section{Homological algebra and stable homotopy theory}
\label{sec:translation}

\subsection*{Background}

Sets (and topological spaces) have a natural diagonal map
$\Delta\co X \to (X^p)^{C_p}$. For an abelian group, we can compose
with the universal multilinear map to get a natural transformation
$A \to (A^{\otimes p})^{C_p}$, given by
$a \mapsto a \otimes \dots \otimes a$. This is a natural
transformation of sets, but not a homomorphism: however, this problem
vanishes modulo the image of the transfer homomorphism
\[
  (A^{\otimes p})_{C_p} \to (A^{\otimes p})^{C_p},
\]
from the $C_p$-quotient to the $C_p$-invariants. It therefore
determines a natural homomorphism from $A$ to the Tate cohomology 
group $\tH^0(C_p;A^{\otimes p})$. This map is the (algebraic)
Tate-valued Frobenius.

To get a chain-level or derived variant of this Frobenius, we would
need to replace the Tate cohomology functor $\tH^0$ by a derived Tate
construction; but now that we are no longer taking a quotient by
transfers, this no longer strictly imposes the homomorphism property.
As a result, in the construction of a derived version of the
Tate-valued Frobenius we will lose the property of staying within
algebra.

Before we introduce the Tate diagonal, we would like to translate the
objects under consideration in \cite{lipshitz-treumann-tate} to stable
homotopy theory.  In this section we will give some brief background
on this translation process. The author claims no originality for the
material in this chapter.

\subsection{Module spectra and chain complexes}

For any ordinary ring $R$, let $\Ch(R)$ be the category of chain
complexes of $R$-modules. For such chain complexes $C$ and $D$, one
can build a function space $\Map_R(C,D)$: start with a set of vertices
given by chain maps $C \to D$, attach paths associated to chain
homotopies, and so on. More concisely, using the Dold--Kan
correspondence one can take the function complex $\Hom_R(C,D)$ and
associate a simplicial set of maps $C \to D$. Because $Ch(R)$ now has
function spaces, one can speak of homotopy limits and colimits in the
category of chain complexes, for example via an associated
$\infty$-category which is made explicit in \cite[\S
1.3.1]{lurie-higheralgebra}. One can also form a localization
$\Ch(R) \to \mathcal{D}(R)$ by inverting the quasi-isomorphisms; on
the level of homotopy categories, this becomes the map from the
classical chain homotopy category of $R$ to the derived category.

Our starting point is the following theorem, which interprets the
category of chain complexes as equivalent to a construction in stable
homotopy theory.

\begin{thm}[{\cite[7.1.1.16, 7.1.2.13]{lurie-higheralgebra}}]
  \label{thm:algebratotopology}
  Let $R$ be a ring. Then there exists an equivalence of
  $\infty$-categories
  \[
    \theta\co \mathcal{D}(R) \too{\sim} \LMod_{HR}
  \]
  between the derived $\infty$-category of differential graded
  $R$-modules and the category of modules over the Eilenberg--Mac Lane
  spectrum $HR$.

  If $R$ is commutative, this extends to an equivalence of symmetric
  monoidal $\infty$-categories, where the source carries the derived
  tensor product $\otimes^L_R$ and the target carries the relative
  smash product $\otimes_{HR}$.
\end{thm}

This result, in several strengths and several guises, has a long
history in the literature and served as a motivation for many
developments. It is present as an analogy in
\cite[5.32]{thomason-ktheory-etale}; as an equivalence between the
derived category of differential graded $R$-modules and the homotopy
category of $HR$-modules in \cite[IV.2.4]{ekmm}; as an equivalence of
model categories in \cite[5.1.6]{schwede-shipley-stablemodules}; and
an extension of this to a monoidal equivalence in
\cite{shipley-dga}. The above formulation is convenient because it
allows us to apply the extensive machinery built in
\cite{lurie-higheralgebra}.

In the description above, the equivalence $\theta$ of
$\infty$-categories preserves structure that can be expressed
in a homotopy-invariant fashion.

\begin{cor}
  Suppose that $K$ is a simplicial set. Then composition with $\theta$
  induces an equivalence of functor $\infty$-categories
  \[
    \Fun(K,\mathcal{D}(R)) \to \Fun(K, \LMod_{HR}).
  \]
\end{cor}

\begin{exam}
  If $K = BG$ is the classifying space of a finite group, maps
  $BG \to \cC$ of $\infty$-categories are coherent actions of $G$ on
  an object of the $\infty$-category $\cC$. This shows that $\theta$
  preserves coherent $G$-actions: chain complexes of $R$-modules with
  a coherent $G$-action are equivalent to $HR$-modules with a coherent
  $G$-action. Chain complexes with strict $G$-action give rise to
  $HR$-modules with coherent $G$-action under $\theta$.
\end{exam}

\begin{cor}
  The functor $\theta$ preserves homotopy limit and colimit
  diagrams.
\end{cor}

\begin{exam}
  Given a chain complex $C$ of $R[G]$-modules, the tensor product
  \[
    E \otimes_{G} C
  \]
  with a projective resolution $E$ of $\mb Z$ over $\mb Z[G]$ is a
  representative for the homotopy colimit in $\mathcal{D}(R)$ of the
  diagram expressing the $G$-action on $C$, which we denote by
  $C_{hG}$. Therefore, it is taken by $\theta$ to a homotopy
  colimit. Similarly, the function complex $\Hom_G(E,C)$ is a
  representative for the homotopy limit $C^{hG}$, and it is taken by
  $\theta$ to a homotopy limit.
\end{exam}

\begin{exam}
  \label{exam:simplicialcoholim}
  Let $f\co \Delta^\op \to \Ch(R)$ represent a simplicial object in
  chain complexes of $R$-modules. Associated to this there is a double
  complex using the standard alternating sign boundary operators, and
  this double complex has an associated totalization. This total
  complex is a representative for the homotopy colimit of the diagram
  $f$. As a result, $\theta$ takes this total complex to a homotopy
  colimit of the diagram $\theta \circ f$.
\end{exam}

\begin{cor}
  If $R$ is commutative and $\oO$ is an $\infty$-operad, $\theta$
  induces an equivalence
  \[
    \Alg_{\oO}(\mathcal{D}(R)) \too{\sim} \Alg_{\oO}(\LMod_{HR})
  \]
  of $\infty$-categories of $\oO$-algebras.
\end{cor}

\begin{exam}
  Suppose $\oO$ is an ordinary operad which is acted on freely by the
  symmetric groups. Then associated to $\oO$ there is an
  $\infty$-operad such that objects with an action of $\oO$ are
  equivalent to algebras over the associated $\infty$-operad. Since
  categories of algebras over $\infty$-operads are invariant under
  symmetric monoidal equivalence, this allows us to translate
  $A_\infty$ and $E_\infty$ algebras between $\Ch(R)$ and
  $\LMod_{HR}$. For example, an associative differential graded
  $R$-algebra $A$ gives rise to an $A_\infty$-algebra $\theta A$ in
  $\LMod_{HR}$. Similarly, differential graded modules and bimodules
  give rise to modules and bimodules over $\theta A$.
\end{exam}

\subsection{Tate constructions}

The classical Tate cohomology of a group $G$ with coefficients in a
module was exported to the category of spectra by Greenlees and May
\cite{greenlees-may-tate} using equivariant stable homotopy
theory. The categories of chain complexes and modules over a ring
spectrum have the special property of being \emph{stable} \cite[\S
1.1.1]{lurie-higheralgebra}: this roughly a higher-categorical lift of
having a triangulated structure. The Tate construction was generalized
to the case of a stable $\infty$-category in
\cite[6.1.6.24]{lurie-higheralgebra}. In this section, we will recall
some of the important properties satisfied by the Tate construction.

\begin{prop}
  Let $G$ be a finite group, $\mathcal{C}$ a stable $\infty$-category
  which admits countable homotopy limits and homotopy colimits, and
  $M$ a $G$-equivariant object of $\mathcal{C}$. Then there is a
  natural \emph{transfer map}\footnote{Other authors
    refer to this as the \emph{norm map}, which we prefer to reserve
    for multiplicative variants.}
  \[
    \Tr\co M_{hG} \to M^{hG}
  \]
  from the derived orbit object to the derived fixed-point
  object. If $M$ is a free object $\oplus_{g \in G} N \simeq
  \prod_{g \in G} N$, this is equivalent to the natural composite
  \[
    \left(\bigoplus_{g \in G} N\right)_{hG} \too{\sim} N
    \too{\sim} \left(\prod_{g \in G} N\right)^{hG}.
  \]
\end{prop}

\begin{defn}
  Let $G$ be a finite group, $\mathcal{C}$ a stable $\infty$-category
  which admits countable limits and colimits, and $M$ a
  $G$-equivariant object of $\mathcal{C}$. We write $M^{tG}$ for the
  cofiber of the transfer $M_{hG} \to M^{hG}$, and refer to it as the
  \emph{$G$-Tate construction on $M$} or simply the Tate construction.
\end{defn}

\begin{prop}
  The Tate construction for a finite group $G$ has the following
  properties.
  \begin{itemize}
  \item It determines a functor $\cC^{BG} \to \cC$, from
    the $\infty$-category of objects of $\cC$ with $G$-action back to the
    $\infty$-category $\cC$.
  \item It preserves finite coproducts, finite products, homotopy
    pushouts, and homotopy pullbacks.
  \item Any functor $\cC \to \dD$ between stable
    $\infty$-categories that preserves countable homotopy limits and
    colimits also preserves Tate constructions. In particular, this is
    true of equivalences.
  \end{itemize}
\end{prop}

A chain complex $M$ of $R[G]$-modules determines an object in $\Ch(R)$
with $G$-action, and as such we can compare the Tate construction
$M^{tG}$ with more classical constructions.

Let $E$ be a projective resolution of $\mb Z$ by finitely generated
free $\mb Z[G]$-modules, and $E^\vee$ the dual complex $\Hom(E,\mb
Z)$. (For instance, we may take $E$ to be the standard bar resolution.)
Then there is a composite
\[
E \to H_0(E) = \mb Z = H_0(E^\vee) \to E^\vee,
\]
where we view $\mb Z$ as a complex concentrated in degree zero, and we
can construct a mapping cone $W$. This complex $W$ is an unbounded
complex of finitely generated free $\mb Z[G]$-modules. Associated to
any chain complex $M$ with $G$-action, there is then a double complex
$W \otimes_G M$ with three realizations.

\begin{defn}
  Let $M$ be a chain complex with $G$-action. We define the following
  chain-level Tate constructions as complexes:
  \begin{align*}
\Tate^\oplus(M)_n &= \bigoplus_{p+q=n} W_p \otimes_G M_q\\
\Tate^\Pi(M)_n &= \prod_{p+q=n} W_p \otimes_G M_q\\
\Tate(M)_n &= \bigcup_N \prod_{p+q = n, p \leq N} W_p \otimes_G M_q
\end{align*}
  The boundary maps in these complexes are the standard boundary maps
  determined by the Leibniz rule $\partial (a \otimes b) = \partial a
  \otimes b + (-1)^{|a|} a \otimes \partial b$.
\end{defn}

\begin{prop}
  The chain-level Tate constructions for the action of $G$ on $M$ have
  the following properties.
  \begin{enumerate}
  \item All three Tate constructions preserve short exact sequences in
    $M$.
  \item There are natural maps $\Tate^\oplus(M) \to \Tate(M) \to
    \Tate^\Pi(M)$.
  \item The map $\Tate^\oplus(M) \to \Tate(M)$ is an isomorphism if
    $M$ is bounded above.
  \item The map $\Tate(M) \to \Tate^\Pi(M)$ is an isomorphism if $M$
    is bounded below.
  \item There is a conditionally convergent Tate cohomology spectral sequence
    \[
      \widehat H^s(G; H_t(M)) \Rightarrow H_{t-s}
      (\Tate(M)).\footnote{Cohomologically minded readers might prefer
        the indexing $\widehat H^p(G; H^q(M)) \Rightarrow H^{p+q}(\Tate(M))$.}
    \]
  \item There is a natural short exact sequence
    \[
      0 \to M_{hG} \to M^{hG} \to \Tate(M) \to 0
    \]
    of complexes, where the first complex is the complex of derived
    coinvariants and the second is the complex of derived invariants.
  \item The object $\Tate(M)$ is a representative for the homotopical
    Tate construction $M^{tG}$. In particular, an equivariant
    quasi-isomorphism $M \to N$ induces an equivalence $\Tate(M) \to
    \Tate(N)$.
  \end{enumerate}
\end{prop}

In particular, for bounded complexes there is no distinction between
these three constructions. However, these three Tate constructions
typically have quite different behavior for unbounded complexes.

\subsection{Algebras and Hochschild complexes}
\label{sec:hochschild}

The equivalence of symmetric monoidal $\infty$-categories between the
category $\mathcal{D}(R)$ and the category $\LMod_{HR}$ allows us to
transport Hochschild complexes because they can be expressed
diagrammatically. Given an associative differential graded $R$-algebra
$A$ with $R$-linear bimodule $M$, the \emph{cyclic bar construction}
is a simplicial chain complex
\[
  Z^R(A,M)\co \Delta^\op \to \Ch(R)
\]
whose associated total complex is the Hochschild complex. If all the
tensor products in this complex are equivalent to the derived tensor
products, then the functor $\theta$ from
Theorem~\ref{thm:algebratotopology} preserves them, and takes this
simplicial diagram to the cyclic bar construction
\[
  Z^{HR}(\theta A, \theta M)\co \Delta^\op \to \LMod_{HR}.
\]
As in Example~\ref{exam:simplicialcoholim}, $\theta$ takes the
associated total complex to the homotopy colimit. The homotopy colimit
is the geometric realization of this cyclic bar construction, which is the
definition of topological Hochschild homology. As a result, we have an
identification:
\[
  \theta (HH^R(A;M)) \simeq \THH^{HR}(\theta A, \theta M).
\]
(If we do not have sufficient flatness of $A$ or $M$ as $R$-modules,
the right-hand side is instead identified with the derived variant of
the Hochschild complex that is used to calculate Shukla homology.)

We would now like to develop the interaction with the cyclic group.
Fix a projective resolution $W$ of $A$ as an $R$-linear
$A$-bimodule. Since any two resolutions are equivalent, there is a
quasi-isomorphism
\[
  B^R(A,A,A) \to W
\]
of $A$-bimodules, where the source is the (total complex associated to
the) two-sided bar resolution relative to $R$. This becomes a
$C_2$-equivariant equivalence
\[
  \left[M \otimes_A B(A,A,A) \otimes_A M\right] \otimes_{A \otimes_R
    A^\op} B(A,A,A) \to \left[M \otimes_A W \otimes_A M\right] \otimes_{A \otimes_R A^\op} W,
\]
where $C_2$ acts by rotational symmetry on the tensor products.
The left-hand side is a bisimplicial object; its homotopy colimit
realizes the Hochschild complex $HH^R(A;M \otimes^L_A M)$. Therefore,
we have an equivalence of Tate constructions
\[
  \Tate(HH^R(A; M \otimes^L_A M)) \simeq \Tate\left(\left[M \otimes_A W
      \otimes_A M\right] \otimes_{A \otimes_R A^\op} W\right)
\]
because $\Tate$ preserves quasi-isomorphisms.

If we now assume that $A$ is homologically smooth, we can choose $W$
to be finitely generated as a complex of bimodules. If $M$ is bounded,
this makes the right-hand complex bounded, and so both become
quasi-isomorphic to the direct-sum Tate complex
\[
  \Tate^\oplus\left(\left[M \otimes_A W \otimes_A M\right] \otimes_{A
    \otimes_R A^\op} W\right).
\]
This last is the complex whose Tate spectral sequence was developed by
Lipshitz and Treumann in \cite{lipshitz-treumann-tate}.

We now apply $\theta$. We find that Lipshitz and Treumann's construction
is carried to a model in stable homotopy theory: the Tate construction
\[
  \left[\THH^{HR}(\theta A;\theta M \otimes_{\theta A} \theta M)\right]^{tC_2}.
\]
This translation now allows us to apply results in stable homotopy
theory to understand Lipshitz and Treumann's Tate spectral sequence.

\section{Envelopes}
\label{sec:envelopes}

We recall that a \emph{multicategory} (also known as a colored operad)
is a category where morphisms may have several inputs: we have a set
$\Hom_\dD(X_1,\dots,X_n;Y)$ of multimorphisms, or alternatively can
describe a map as being of the form $\{X_s\}_{s \in S} \to Y$ with $S$
a finite index set. Every symmetric monoidal category $\dD$ has an
underlying multicategory $U\dD$: we define maps
$\{X_s\}_{s \in S} \to Y$ to be the same as maps
$\bigotimes^\dD X_s \to Y$. In the other direction, associated to a
multicategory $\cC$ there is a symmetric monoidal category $\Env(\cC)$
called the \emph{symmetric monoidal envelope}.
\begin{itemize}
\item The objects of $\Env(\cC)$ are formal tuples
  $(S,\{X_s\}_{s \in S})$ of a finite set and an $S$-indexed set of
  objects of $\cC$, representing a formal tensor $\bigotimes^\cC
  X_s$.
\item The morphisms $(S,\{X_s\}) \to (T, \{Y_t\})$ in $\Env(\cC)$ are
  pairs of a map $f\co S \to T$ and a collection of maps
  $g_t\co \{X_s\}_{s \in p^{-1}(t)} \to Y_t$ in the multicategory
  $\cC$.
\end{itemize}
By construction, there is an equivalence between multifunctors
$\cC \to U\dD$ and symmetric monoidal functors $\Env(\cC) \to \dD$;
there is also a forgetful functor from $\Env(\cC)$ to the category of
finite sets.

The coherent version of this construction is described in
\cite[\S2.2.4]{lurie-higheralgebra}. The analogues of multicategories
are \emph{$\infty$-operads}, and an $\infty$-operad $\cC^\otimes$
has an associated symmetric monoidal envelope $\Env(\cC)$. The
universal property of the symmetric monoidal envelope is
\cite[2.2.4.9]{lurie-higheralgebra}: for any symmetric monoidal
$\infty$-category $\dD$, there is an equivalence between symmetric
monoidal functors $\Env(\cC) \to \dD$ and maps of $\infty$-operads
$\cC^\otimes \to \dD^\otimes$.

Here are some important properties of this construction.
\begin{itemize}
\item There is a natural symmetric monoidal functor
  $p\co \Env(\cC) \to \Fin$ to the category of finite sets, whose
  fiber over $S$ is equivalent to $\prod_{s \in S} \cC$, and under
  this correspondence the space of maps
  $\{X_s\}_{s \in S} \to \{Y_t\}_{t \in T}$ over a given map
  $f\co S \to T$ is equivalent to a product
  $\prod_{t \in T} \Map_{\cC^\otimes}(\{X_s\}_{s \in f^{-1} (t)}, Y_t)$ of
  mapping spaces.
\item The adjunction gives every symmetric monoidal $\infty$-category
  $\dD$ a natural symmetric monoidal functor
  $\bigotimes^\dD\co \Env(\dD) \to \dD$, sending $\{X_s\}_{s \in S}$
  to $\bigotimes^\dD_{s \in S} X_s$.
\end{itemize}

\begin{exam}
  For any associative algebra $A$ in a symmetric monoidal
  $\infty$-category $\cC$ with right module $M$ and left module $N$,
  the two-sided bar construction $B(M,A,N)$ can be lifted from a
  simplicial object in $\cC$ to a simplicial object in
  $\Env(\cC)$. This is precisely because the face and degeneracy
  operations involve only multiplication, permutation of factors, and
  insertion of units, and never any nontrivial maps of the form
  $X \to Y_1 \otimes \dots \otimes Y_n$.
\end{exam}

\begin{exam}
  Similarly, given an associative algebra $A$ in a symmetric monoidal
  $\infty$-category $\cC$, there exists a lift of the cyclic bar
  construction $Z(A)$ from a simplical object in $\cC$ to a simplicial
  object in $\Env(\cC)$. For an algebra $A$ with a bimodule $M$, the
  same is true for $Z(A,M)$, the cyclic bar construction with
  coefficients.
\end{exam}

\section{Pushforward}

If $\cC$ is an ordinary symmetric monoidal category, the counit of the
adjunction between symmetric monoidal categories and multicategories
takes the form of a tensor functor
\[
  \bigotimes^\cC\co \Env(U\cC) \to \cC,
\]
taking an $S$-indexed tuple $\{X_s\}_{s \in S}$ to
$\bigotimes^\cC_{s \in S} X_s$. More generally, this admits a relative
version of the following type. Suppose that we have an $S$-indexed
tuple $X = \{X_s\}_{s \in S}$ of objects of $\cC$ and a map
$\varphi\co S \to T$ of finite sets. Then, associated to this, there
is a $T$-indexed tuple
\[
  \varphi_!(X) = \left\{\bigotimes^\cC_{s \in f^{-1}(t)}
    X_s\right\}_{t \in T}
\]
(which we call the \emph{fiberwise tensor}) together with a map
$X \to \varphi_!(X)$ over $\varphi$.

Each $\varphi_!$ is a functor $\cC^S \to \cC^T$. However, there is
further structure: functoriality extends to categories of lifts. Given
any category $I$ and a functor $f\co I \to \Fin$ from $I$ to the
category of finite sets, there is a category $\Env(\cC,f)$ of
\emph{lifts of $f$}: functors $F\co I \to \Env(\cC)$ lifting $f$, and
therefore associating to each $i \in I$ an $f(i)$-indexed tuple
$\{F(i)_s\}_{s \in f(i)}$. Any natural transformation
$\varphi\co f \to g$ of diagrams of sets gives rise, for each
$i \in I$, to a functor $\varphi(i)_!$ from $f(i)$-indexed tuples to
$g(i)$-indexed tuples, and these assemble together into a global
fiberwise tensor functor $\varphi_!\co \Env(\cC, f) \to \Env(\cC,
g)$. Up to coherent natural isomorphism, this construction preserves
composition in $\varphi$ and is equivalent to a functor
\[
  \Env(\cC,-)\co \Fun(I,\Fin) \to \Cat.
\]
Even more, the categories $\Env(\cC,f)$ are symmetric monoidal under
the pointwise tensor product and the functors $\varphi_!$ are
symmetric monoidal.

In this section we will discuss the construction of the analogous
fiberwise tensor in the context of $\infty$-categories. The main
result is Proposition~\ref{prop:functorialfiberwise}, which is an
``unstraightened'' variant of the symmetric monoidal functoriality of
$\Env(\cC,-)$. The results of this section make heavier use of
technical details from \cite{lurie-higheralgebra}, and the reader
without close familiarity may not have need for these details.

We first require some general intermediate results.

\begin{lem}
  \label{lem:cocartenvelope}
  Suppose that $\cC^\otimes \to \dD^\otimes$ is a coCartesian
  fibration of $\infty$-operads. Then the functor $\Env(\cC)^\otimes
  \to \Env(\dD)^\otimes$ is a coCartesian fibration of symmetric
  monoidal $\infty$-categories.
\end{lem}

\begin{proof}
  The symmetric monoidal envelope $\Env(\cC)^\otimes$ is defined in
  \cite[2.2.4.1]{lurie-higheralgebra} as the fiber product
  \[
    \cC^\otimes \times_{\Fin_*} \mathrm{Act}(\Fin_*)
  \]
  with the subcategory spanned by active morphisms. The result follows
  because fibration conditions are stable under base-change.
\end{proof}

\begin{lem}
  \label{lem:cocartfunctor}
  Suppose that $\cC^\otimes \to \dD^\otimes$ is a coCartesian
  fibration of symmetric monoidal $\infty$-categories and that
  $\oO^\otimes$ is an $\infty$-operad. Then the functor
  \[
    \Alg_\oO(\cC) \to \Alg_\oO(\dD)
  \]
  lifts to a symmetric monoidal coCartesian fibration under the
  pointwise tensor product of $\oO$-algebras from
  \cite[3.2.4.4]{lurie-higheralgebra}.
\end{lem}

\begin{proof}
  For $\dD^\otimes$ a symmetric monoidal $\infty$-category and
  $\oO^\otimes$ an $\infty$-operad, $\Alg_\oO(\dD)$ is the full
  subcategory of $\Fun(\oO^\otimes, \dD^\otimes)$ spanned by the maps
  of $\infty$-operads. The definition of the pointwise tensor product
  in \cite[3.2.4.4]{lurie-higheralgebra} is as follows: maps
  $K \to \Alg_\oO(\dD)^\otimes$ over a fixed map $K \to \Fin_*$ are
  equivalent to commutative diagrams
  \[
    \xymatrix{
      K \times \oO^\otimes \ar[r] \ar[d] &
      \dD^\otimes \ar[d] \\
      \Fin_* \times \Fin_* \ar[r]^-\wedge &
      \Fin_*
    }
  \]
  such that the top map restricts to a map of $\infty$-operads for any
  vertex of $K$.

  By adjunction, then, the identity self-functor of
  $\Alg_\oO(\dD)^\otimes$ determines a commutative diagram
  \[
    \xymatrix{
      \Alg_\oO(\dD)^\otimes \times \oO^\otimes \ar[r] \ar[d] &
      \dD^\otimes \ar[d] \\
      \Fin_* \times \Fin_* \ar[r]^-\wedge &
      \Fin_*.
    }
  \]
  The topmost map sends pairs of inert morphisms in
  $\Alg_\oO(\dD)^\otimes \times \oO^\otimes$ to inert morphisms in
  $\dD^\otimes$ \cite[3.2.4.3, (2)]{lurie-higheralgebra}, and thus
  it is a \emph{bifunctor} of $\infty$-operads in the sense of
  \cite[2.2.5.3]{lurie-higheralgebra}.

  We now apply \cite[3.2.4.3]{lurie-higheralgebra} to the bifunctor
  $\Alg_\oO(\dD)^\otimes \times \oO^\otimes \to \dD^\otimes$ and the
  coCartesian fibration $\cC^\otimes \to \dD^\otimes$. This shows that,
  under the definition from \cite[3.2.4.1]{lurie-higheralgebra},
  there is a coCartesian fibration
  \[
    \Alg_{\oO/\dD}(\cC)^\otimes \to \Alg_\oO(\dD)^\otimes.
  \]
  However, unravelling the definition of the source we find that this
  is the natural functor
  \[
    \Alg_{\oO}(\cC)^\otimes \to \Alg_\oO(\dD)^\otimes
  \]
  under the pointwise monoidal structure. In particular, the fiber
  over a map of $\infty$-operads $f\co \oO^\otimes \to \dD^\otimes$ is
  the $\infty$-category of lifts $\oO^\otimes \to \cC^\otimes$.

  Moreover, by \cite[3.2.4.3, (4)]{lurie-higheralgebra}, a morphism
  $\alpha\co A \to B$ in $\Alg_\oO(\cC)^\otimes$ is coCartesian if and
  only if, for any $X \in \oO$, the natural transformation the map
  $A(X) \to B(X)$ of $\cC^\otimes$ is a coCartesian lift of
  its image in $\dD^\otimes$.
\end{proof}

% \begin{rmk}
%   The symmetric monoidal structure on $\Fin$ is the coproduct of
%   finite sets, and the induced symmetric monoidal structure on
%   $\Fun(K,\Fin)$ is the coproduct of diagrams of finite sets. Every
%   functor $f\co K \to \Fin$ is naturally a commutative monoid under
%   this structure, and because $\Env(\cC,-)$ is defined by taking
%   categories of sections it takes coproducts to products: it is a symmetric
%   monoidal functor from $\Fun(K,\Fin)$ to symmetric monoidal
%   $\infty$-categories.

%   More explicitly, the symmetric monoidal structure on $\Env(\cC,f)$
%   comes from the map $\Env(\cC) \times_{\Fin} \Env(\cC) \to \Env(\cC)$
%   given by the tensor product in $\cC$ fiberwise.
% \end{rmk}

\begin{prop}
  \label{prop:functorialfiberwise}
  If $\cC$ is a symmetric monoidal $\infty$-category, the map of
  functor categories $\Fun(K,\Env(\cC)) \to \Fun(K,\Fin)$ extends, up
  to equivalence, to a symmetric monoidal coCartesian fibration.
\end{prop}

\begin{proof}
  Because $\cC$ is symmetric monoidal, we have a coCartesian fibration
  $\cC^\otimes \to \Fin_*$ and hence a coCartesian fibration
  \[
    \Env(\cC)^\otimes \to \Env(\Fin)^\otimes = \Fin^\amalg
  \]
  by Lemma~\ref{lem:cocartenvelope}.
  
  Given a simplicial set $K$, viewed as a simplicial set over $\Fin_*$
  via $K \to \{1\} \subset \Fin_*$, let
  $K \to \mathcal{K}^\otimes \to \Fin_*$ be a fibrant replacement in
  the $\infty$-operadic model structure
  \cite[2.1.4.6]{lurie-higheralgebra}; $\mathcal{K}^\otimes$ is an
  $\infty$-operad. This has the property that for any $\infty$-operad
  $\cC^\otimes$, restricting maps of $\infty$-operads
  $\mathcal{K}^\otimes \to \dD^\otimes$ to functors
  $K \to \dD$ gives an equivalence of functor categories
  \[
    \Alg_{\mathcal{K}}(\dD) \simeq \Fun(K,\dD).
  \]
  For a symmetric monoidal $\infty$-category $\cC$, we then get a
  commutative diagram
  \[
    \xymatrix{
      \Alg_{\mathcal{K}}(\Env(\cC)) \ar[r] \ar[d] &
      \Fun(K,\Env(\cC)) \ar[d] \\
      \Alg_{\mathcal{K}}(\Fin) \ar[r] &
      \Fun(K,\Fin).
    }
  \]
  The horizontal maps are equivalences, and the left-hand vertical map
  extends to a symmetric monoidal coCartesian fibration by
  Lemma~\ref{lem:cocartfunctor}.
\end{proof}

\begin{rmk}
  The straightening of this coCartesian fibration is a symmetric
  monoidal functor
  \[
    \Env(\cC,-)\co \Fun(K,\Fin) \to \Cat_\infty.
  \]
  This sends $f\co K \to \Fin$, up to equivalence, to the category
  $\Env(\cC,f)$ of lifts $F\co K \to \Env(\cC)$ of $f$ and sends a
  natural transformation $\varphi\co f \to g$ to a functor
  $\varphi_!\co \Env(\cC,f) \to \Env(\cC,g)$. Functoriality says that
  this respects composition in $\varphi$, and the coCartesian property
  means that every natural transformation $F \to G$ over $\varphi$
  factors essentially uniquely through a functor $\varphi_! F \to G$
  over $\varphi$. Moreover, the description of coCartesian morphisms
  leads to the following: the map $\alpha\co F \to \varphi_! F$ is
  characterized by the property that, for any object $k$ of $K$, the map
  $\alpha(k)\co F(k) \to (\varphi_! F)(k)$ of $\cC^\otimes$ is a
  coCartesian lift of the underlying map $\varphi(k)\co f(k) \to
  g(k)$. In other words, if $F(k) = \{X_s\}_{s \in f(k)}$, then
  \[
    (\varphi_! F)(k) =
    \left\{\bigotimes^\cC_{\varphi(s) = t} X_s\right\}_{t \in g(k)}.
  \]

  % \todo[inline]{mess}
  % bifunctor $\oO^\otimes \times {\oO'}^\otimes \to {\oO''}^\otimes$
  % and a fibration $q\co \cC^\otimes \to {\oO''}^\otimes$

  % FOR US: $\Alg_\oO(\dD)^\otimes \times \oO^\otimes \to \dD^\otimes$
  % and $q: \cC^\otimes \to \dD^\otimes$

  % induced $p\co \Alg_{\oO'/\oO''}(\cC)^\otimes \to \oO^\otimes$

  % FOR US: $p: \Alg_{\oO / \dD}(\cC)^\otimes \to \Alg_{\oO}(\dD)^\otimes$
  
  % a morphism $\alpha$ of $\Alg_{\oO'/\oO''}(\cC)^\otimes$ is
  % $p$-coCartesian iff, for any $X \in \oO'$, the image $\alpha(X)$ is
  % $q$-coCartesian.

  % FOR US: a morphism $\alpha \in \Alg_{\oO/\dD}(\cC)^\otimes$ is
  % $p$-coCart iff image of any $X$ is coCart
  % \todo[inline]{end mess}
  
  The unit of $\Env(\dD,f)$ is a functor
  $\mb I_f\co K \to \Env(\dD)$, sending an object $k$ to an indexed
  tuple $\{\mb I_\dD\}_{s \in f(k)}$ of copies of the unit for the
  tensor product.
  % \todo{check this wtf why do I care}This respects restrictions in the source
  % category as follows. Given a natural transformation
  % $\varphi\co f \to g$ and a functor $q\co L \to K$, the natural map
  % \[
  %   (q^* \varphi)_! (q^* F) \to
  %   q^*(\varphi_! F)
  % \]
  % is an equivalence of functors $L \to \cC$ for any $F$.
\end{rmk}

% \begin{cor}
%   Given a natural transformation $\theta$ of functors
%   $f,g\co K \to \Fin$, the fiberwise tensor extends to a natural
%   transformation $F \to \theta_! F$ over $\theta$. Any natural
%   transformation $F \to G$ over $\theta$ factors through it
%   essentially uniquely. In particular, there is always a natural
%   equivalence $\theta_! \varphi_! = (\theta \varphi)_!$.
% \end{cor}

\begin{exam}
  \label{exam:pushcyclic}
  The quotient map $\varphi\co \Delta^1 \to S^1$ of simplicial sets is
  a natural transformation of functors $\Delta^\op \to \Fin$. For any
  algebra $A$ in $\cC$ with right module $M$ and left module $N$, the
  fiberwise tensor $\varphi_! B(M,A,N)$ of the two-sided bar
  construction is the cyclic bar construction $Z(A, N \otimes M)$
  associated to the $A$-bimodule $N \otimes M$.
\end{exam}

\begin{exam}
  \label{exam:pushforward}
  For any $K$, there is a constant functor $\ast\co K \to \Fin$ whose
  value is a singleton, and lifts of this to $\Env(\cC)$ are
  equivalent to functors $K \to \cC$. There always a natural
  transformation $\tau$ from any functor $f\co K \to \Fin$ to the
  constant functor with value $*$, and under these identifications the
  functor $\tau_! F$ is the tensor product functor
  $\bigotimes^\cC \circ F\co K \to \cC$. In particular, functoriality
  of fiberwise tensor tells us that we have natural equivalences
  \[
    \bigotimes^\cC \circ \varphi_! F = \bigotimes^\cC \circ F.
  \]
\end{exam}

\section{Realization}

We now consider homotopy colimits. For an $\infty$-category $\cC$, a
given diagram $K \to \cC$ may have a homotopy colimit. If $\cC$ has
all $K$-indexed homotopy colimits these can be made into a functorial
homotopy colimit $\Fun(K,\cC) \to \cC$. If $\cC$ is also symmetric
monoidal, we can also construct natural transformations
\[
  \hocolim_i (F(i) \otimes G(i)) \to \hocolim_{i,j} (F(i) \otimes
  G(j)) \to (\hocolim_i F(i)) \otimes (\hocolim_j G(j));
\]
but these maps are not natural equivalences without making further
assumptions. The first map is an equivalence if the index category is
sifted (the diagonal $\Delta\co K \to K \times K$ is cofinal
\cite[5.5.8.1]{lurie-htt}), and the second map is an equivalence if
the monoidal product of $\cC$ preserves homotopy colimits in each
variable separately.

\begin{prop}
  \label{prop:laxhocolim}
  Suppose that $K$ is sifted, that $\cC$ is a symmetric monoidal
  $\infty$-category with $K$-indexed colimits, and that the symmetric
  monoidal structure preserves $K$-indexed colimits in each
  variable. Then there is a functor
  \[
    \hocolim_K\co \Fun(K,\cC) \to \cC
  \]
  that is strong symmetric monoidal.
\end{prop}

A commutative algebra object $A$ in an $\infty$-operad determines a
symmetric monoidal functor $A\co \Fin \to \Env(\cC)$, where the image
of a finite set $S$ has a chosen equivalence with a constant indexed
tuple $\{A\}_{s \in S}$.
\begin{defn}
  Suppose that $K$ is sifted, that $\cC$ is a symmetric monoidal
  $\infty$-category with $K$-indexed colimits, and that the symmetric
  monoidal structure preserves $K$-indexed colimits in each
  variable. Let $A$ be a commutative algebra object in $\cC$. We
  define
  \[
    A^{\otimes X} = \hocolim_K \left(\bigotimes^\cC \circ A \circ
      X\right)
  \]
  for any functor $X\co K \to \Fin$.
\end{defn}

\begin{exam}
  Suppose that $X$ is a simplicial finite set, viewed as a functor
  $X\co \Delta^\op \to \Fin$. Then $A^{\otimes X}$ can be identified
  with the Loday construction $A^{\otimes |X|}$.
\end{exam}

\begin{prop}
  Suppose that $K$ is sifted, that $\cC$ is a symmetric monoidal
  $\infty$-category with $K$-indexed colimits, and that the symmetric
  monoidal structure preserves $K$-indexed colimits in each
  variable. Let $\varphi\co f \to g$ be a natural transformation of
  functors $K \to \Fin$ and $F\co K \to \Env(\cC)$ a lift of
  $f$. Then there is a natural equivalence
  \[
    \hocolim_K \left(\bigotimes^\cC \circ F\right) \to \hocolim_K
    \left(\bigotimes^\cC \circ \varphi_! F\right).
  \]
\end{prop}

\begin{proof}
  This follows from the equivalence between $\bigotimes^\cC \circ
  \varphi_! F$ and $\bigotimes^\cC \circ F$ from
  Example~\ref{exam:pushforward}.
\end{proof}

\begin{rmk}
  For convenience, in the remainder of the paper we will typically
  make the stronger assumption that the category in question is
  symmetric monoidal presentable.
\end{rmk}

\begin{exam}
  If $M$ is a right $A$-module and $N$ is a left $A$-module, then the
  identification between the pushforward of the two-sided bar
  construction and a cyclic bar construction from
  Example~\ref{exam:pushcyclic} gives us an equivalence
  \[
    M \otimes_A N \simeq \THH(A;N \otimes M).
  \]
\end{exam}

\section{Adjoints and algebras}

Envelopes are functorial: for a map of $\infty$-operads
$R\co \dD^\otimes \to \cC^\otimes$, there is an induced symmetric
monoidal functor $\Env(R)\co \Env(\dD) \to \Env(\cC)$. Further, if
$\cC$ and $\dD$ are symmetric monoidal $\infty$-categories, a map of
$\infty$-operads $\dD^\otimes \to \cC^\otimes$ encodes a lax symmetric
monoidal functor $R\co \dD \to \cC$. There is a resulting natural
transformation
\[
  \bigotimes^\cC \circ \Env(R) \to R \circ \bigotimes^{\dD},
\]
and the functor $R$ is strong symmetric monoidal precisely when this
is a natural equivalence. Moreover, the right adjoint to a
strong symmetric monoidal functor is lax symmetric monoidal
\cite[7.3.2.7]{lurie-higheralgebra}. When we apply this to categories
of lifts, we find the following result.

\begin{prop}
  \label{prop:laxstruct}
  Suppose that $R\co \dD \to \cC$ is a lax symmetric monoidal
  functor, and $f\co K \to \Fin$ is a fixed functor. Then there are
  induced lax symmetric  monoidal functors $\Env(R,f)\co \Env(\dD,f) \to
  \Env(\cC,f)$ and $R\co \Fun(K,\dD) \to \Fun(K,\cC)$, together with a
  lax symmetric monoidal natural transformation
  \[
    \bigotimes^\dD \circ \Env(R) \circ F \to R \circ \bigotimes^\cC
    \circ F
  \]
  for $F \in \Env(\dD,f)$.
\end{prop}

In general, a lax symmetric monoidal functor $F$ takes a commutative
algebra $A$ to a commutative algebra $F(A)$ and an $A$-module $M$ to
an $F(A)$-module $F(M)$. A lax symmetric monoidal natural
transformation $F \to G$ takes a commutative algebra $A$ to a natural
map of commutative algebras $F(A) \to G(A)$ and an $A$-module $M$ to a
map $F(M) \to G(M)$ of $F(A)$-modules. In the case of these envelope
categories, this specializes in the following way.

\begin{prop}
  Suppose that we have a lax symmetric monoidal functor
  $R\co \dD \to \cC$.  Given $f\co K \to \Fin$, then there is a
  natural map of commutative algebras
  \[
    \bigotimes^\cC R(\mb I_f) \to R(\mb I_{\dD})
  \]
  in $\Fun(K,\cC)$. The functor $\bigotimes^\cC \circ R\co \Env(\dD,f)
  \to \Fun(K,\cC)$ lifts to the category of $\bigotimes^\cC R(\mb
  I_f)$-modules; the functor $R \circ \bigotimes^\dD\co \Env(\dD,f)
  \to \Fun(K,\cC)$ lifts to the category of $R(\mb
  I_\dD)$-modules; and the natural transformation of
  Proposition~\ref{prop:laxstruct} lifts to natural a map of modules.
\end{prop}

We now compose this with the natural transformation $\hocolim_K
\circ R \to R \circ \hocolim_K$.

\begin{prop}
  \label{prop:generallax}
  Suppose that $K$ is sifted, $f\co K \to \Fin$ is fixed, and that
  $R\co \dD \to \cC$ is a lax symmetric monoidal functor between
  symmetric monoidal presentable $\infty$-categories. Then there is a
  natural map of commutative algebras
  \[
    R(\mb I_f)^{\otimes f} \to R(\mb I_{\dD})
  \]
  in $\cC$. The functor
  $\hocolim_K \circ \bigotimes^\cC \circ R\co \Env(\dD,f) \to \cC$
  lifts to the category of $R(\mb I_\dD)^{\otimes f}$--modules; the
  functor
  $R \circ \hocolim_K \circ \bigotimes^\dD\co \Env(\dD,f) \to \cC$
  lifts to the category of $R(\mb I_\dD)$-modules; there is an induced
  transformation
  \[
    \hocolim_K \left(\bigotimes^\cC \circ R \circ F\right) \to 
    R \hocolim_K \left(\bigotimes^\dD\circ F\right)
  \]
  of $R(\mb I_\dD)^{\otimes f}$-modules.
\end{prop}

\begin{defn}
  Suppose that $K$ is sifted, $f\co K \to \Fin$ is fixed, and that
  $R\co \dD \to \cC$ is a lax symmetric monoidal functor between
  symmetric monoidal presentable $\infty$-categories. The
  \emph{base-change map} is the natural transformation
  \[
    R(\mb I_\dD) \mathop\otimes_{R(\mb I_\dD)^{\otimes f}}
    \left(\hocolim_K \bigotimes^{\cC} \circ R \circ F\right) \to
    R\left(\hocolim_K \bigotimes^{\dD} \circ F\right)
  \]
  of functors $\Env(\dD,f) \to \LMod_{R(\mb I_{\dD})}$, adjoint to the
  map of Proposition~\ref{prop:generallax}.
\end{defn}

We now specialize this to the case where $\dD$ is the category
$\LMod_A$ of left
modules over a fixed commutative algebra object.

\begin{thm}
  Let $A$ be a commutative algebra in a symmetric monoidal presentable
  $\infty$-category $\cC$. Suppose that $K$ is sifted,
  $f\co K \to \Fin$ is fixed, and that $R\co \LMod_A \to \cC$ is the
  forgetful functor. Then $R$ is lax symmetric monoidal, and the
  base-change map is a natural equivalence.
\end{thm}

\begin{proof}
  The functor $R$ is right adjoint to the strong symmetric monoidal
  functor $X \mapsto A \otimes X$, and hence is lax monoidal
  \cite[7.3.2.7]{lurie-higheralgebra}. The functor $R$ also preserves
  homotopy limits and colimits \cite[4.2.3.3,
  4.2.3.5]{lurie-higheralgebra}.

  We will first prove that the base-change map is an equivalence in
  the case where $K = *$ is the trivial category. In this case,
  without loss of generality, the map $f$ is a choice of a finite set
  $S$ and a lift $F$ is equivalent to an $S$-indexed tuple $\{M_s\}$
  of left $A$-modules. The base-change map is the map
  \[
    A \mathop\otimes_{A^{\otimes S}}
    \left(\bigotimes M_s\right) \to {\bigotimes}^A M_s.
  \]
  The base-change map is an equivalence whenever each $M_s$ is an extended
  module of the form $A \otimes X_s$ for some $X_s$.
  
  The natural augmentation of left $A$-modules $B(A,A,M) \to M$ from
  the two-sided bar construction gives rise to a diagram 
  \[
    \xymatrix{
      \hocolim_{\Delta^\op} A \otimes_{A^{\otimes S}} \bigotimes
      B(A,A,M_s) \ar[r] \ar[d] &
      \hocolim_{\Delta^\op} \bigotimes^A B(A, A, M_s) \ar[d] \\
      A \otimes_{A^{\otimes S}} \bigotimes M_s \ar[r] &
      \bigotimes^A M_s.
    }
  \]
  The top map is a homotopy colimit of a diagram of equivalences
  because the bar construction levelwise consists of extended modules.
  Since $\cC$ and $\LMod_A$ are presentable symmetric monoidal, by
  definition the tensor product preserves homotopy colimits in each
  variable and sifted homotopy colimits in general; therefore, the
  left and right maps are equivalences. The bottom map is then
  an equivalence.

  Now suppose that $K$ is a general sifted index category with map
  $f\co K \to \Fin$. The base-change map is a map
  \[
    A \mathop\otimes_{\hocolim_{k \in K} A^{\otimes f(k)}}
    \left(\hocolim_{k \in K}\bigotimes_{s \in f(k)} F(k)_s\right)
    \to
    \hocolim_{k \in K}
    \bigotimes^A_{k \in K,s \in f(k)} F(k)_s
  \]
  Since $K$ is sifted, and the forgetful functor preserves homotopy
  colimits, we can rewrite both sides as homotopy colimits indexed by
  $k \in K$. The base-change map is then equivalent to the homotopy
  colimit of the base-change maps indexed by $f(k)$, which we already
  showed to be equivalences.
\end{proof}
%\todo{OMG examples}

\section{Shape}

\begin{defn}
  Let $\cC^\otimes$ be an $\infty$-operad with symmetric monoidal
  envelope $\Env(\cC)$. Given a functor $X\co K \to \Env(\cC)$, the
  \emph{shape} of $X$, denoted by $|X|$, is the homotopy colimit of
  the composite functor to the category $\Spaces$ of spaces:
  \[
    K \to \Env(\cC) \to \Fin \subset \Spaces
  \]
\end{defn}

\begin{exam}
  Suppose that $X\co \Delta^\op \to \Env(\cC)$ is a simplicial
  object. Then the composite $\Delta^\op \to \Env(\cC) \to \Fin$ is a
  simplicial finite set, which can be identified with the shape
  $|X|$.
\end{exam}

\begin{exam}
  Suppose $A$ is an associative algebra in $\cC$. Then the functor
  $\Env(\cC) \to \Fin$ takes $A$ to the associative algebra $\ast$
  under coproduct. The cyclic bar construction $Z^{\otimes}(A)$ maps
  to the cyclic bar construction $Z^{\amalg} (\ast)$ and the
  associated shape is $S^1$.

  Similarly, suppose $A$ is an associative algebra with a left module
  $N$ and a right module $M$. Then $\Env(\cC) \to \Fin$ takes the
  two-sided bar construction $\TBar^\otimes(M,A,N)$, whose homotopy
  colimit is $M \otimes_A N$, to the two-sided bar construction
  $\TBar^{\amalg}(\ast, \ast, \ast)$, which is isomorphic to the
  standard simplex $\Delta^1$.
\end{exam}

\begin{defn}
  Let $\cC^\otimes$ be an $\infty$-operad and $Y$ be a Kan complex. We
  define
  \[
    \Env(\cC)_{/Y} = \Env(\cC) \times_{\Spaces} \Spaces_{/Y}.
  \]
\end{defn}

\begin{prop}
  Given an $\infty$-operad $\cC$ and a functor $X\co K \to \Env(\cC)$,
  the category of lifts of $X$ to a functor $\wt X\co K \to
  \Env(\cC)_{/Y}$ is equivalent to the space of maps $f\co |X| \to Y$.
\end{prop}

\begin{proof}
  By definition of the fiber product, lifts $\wt X$ are
  equivalent to lifts of the composite $K \to \Spaces$ to
  $\Spaces_{/Y}$; by definition of the slice category, these are
  equivalent to lifts of $K \to \Spaces$ to natural maps from the
  diagram $K$ to $Y$. However, the universal property of homotopy
  colimits precisely asserts that these extensions are equivalent to
  maps $|X| \to Y$.
\end{proof}

\section{Free $G$-sets}

In this section we will fix a finite group $G$ and let $BG$ be a
Kan complex classifying principal $G$-bundles.

\begin{defn}
  Let $Fin$ be the category of finite sets, and $\Free(G)$ the
  category of finite free left $G$-sets and equivariant maps.
  % There is a \emph{quotient} functor $q\co \Free(G) \to Fin$ that
  % sends a free $G$-set $S$ to the quotient $\bar S = G \setminus S$,
  % and an \emph{underlying} functor $u\co \Free(G) \to Fin$ that
  % forgets the $G$-action.
\end{defn}

\begin{rmk}
  Both categories are symmetric monoidal under disjoint union, but
  their symmetric monoidal structures are particularly simple.
  
  The category $Fin$ is the symmetric monoidal envelope of the
  terminal multicategory $\{*\}$. In particular, any multicategory $C$
  has a canonical symmetric monoidal functor $\Env(C) \to Fin$,
  sending $\{x_s\}_{s \in S}$ to the indexing set $S$.  Moreover, the
  one-point set $*$ is an algebra in $Fin$ and as such is classified
  by a symmetric monoidal functor $\Env(Assoc) \to Fin$.

  Similarly, the category $\Free(G)$ is the symmetric monoidal
  envelope of a one-object multicategory with underlying category
  $BG$.
\end{rmk}

\begin{prop}
  Let $\Spaces$ be the category of spaces, with $Fin$ viewed as a full
  subcategory. The functor $\Free(G) \to Fin_{/BG}$, given by $X
  \mapsto (p_X\co EG \times_G X \to BG)$, induces an equivalence of
  $\infty$-categories
  \[
    \Free(G) \to Fin_{/BG} = Fin \times_{\Spaces} \Spaces_{/BG}.
  \]
  In particular, the space of lifts of a functor $X\co I \to Fin$
  to a functor $\wt X\co I \to \Free(G)$ is equivalent to the
  space of maps $f\co \hocolim_I X \to BG$, classifying principal
  $G$-bundles on the homotopy colimit.
\end{prop}

\begin{proof}
  Because $BG$ is path-connected, an object $S \to BG$ is equivalent in
  $Fin_{/BG}$ to the image of $G \times S$. Therefore, this functor is
  essentially surjective, and so it suffices to show that it is fully
  faithful. This amounts to the assertion that for finite free
  $G$-sets $X$ and $Y$, the diagram
  \[
    \xymatrix{
      \Map_G(X,Y) \ar[r] \ar[d] &
      \Map(EG \times_G X, EG \times_G Y) \ar[d] \\
      \{p_X\} \ar[r] & \Map(EG \times_G X, BG)
    }
  \]
  is a homotopy pullback diagram.

  This diagram decomposes as a product diagram over the orbits of $X$,
  and so it suffices to take $X = G$. However, in this case we recover
  something equivalent to the standard homotopy pullback diagram
  \[
    \xymatrix{
      Y \ar[r] \ar[d] &
      EG \times_G Y \ar[d] \\
      \ast \ar[r] & BG.\qedhere
    }
  \]
\end{proof}

\begin{cor}
  There is an equivalence
  \[
    \Free(G) \times_{Fin} \Env(\cC) \simeq \Env(\cC)_{/BG}.
  \]
\end{cor}

As a result, we write diagrams $K \to \Env(\cC)_{/BG}$ as pairs
$(X,f)$ of a functor $X\co K \to \Env(\cC)$ and a classifying map
$f\co |X| \to BG$.

\section{Unwinding}

\begin{prop}
  For an $\infty$-operad $\cC$, there is a symmetric monoidal
  \emph{fiberwise tensor power} functor
  \[
    \psi\co \Env(\cC)_{/BG} \to \Env(\cC)^{BG}.
  \]
\end{prop}

Informally, the functor $\psi$ sends a free $G$-set $S$ and an
$\bar S$-indexed family $\{c_{\bar s}\}_{\bar s \in \bar S}$ to the
$S$-indexed family $\{c_s\}_{s \in S}$ with its $G$-action.

\begin{proof}
  This is stated for $G = C_p$ in \cite[III.3.6]{scholze-nikolaus-tc},
  but the proof does not make use of any structure particular to this
  group. We briefly recall their method.

  The category $\Env(\cC)$ is symmetric monoidal, and the category
  $\Fun_\otimes(\Env(\cC), \Env(\cC))$ of symmetric monoidal functors
  inherits a pointwise symmetric monoidal structure. The inclusion of
  the identity functor $id$ induces a symmetric monoidal functor from
  $\Fin$, the free symmetric monoidal $\infty$-category on $\{id\}$, to
  $\Fun_\otimes(\Env(\cC), \Env(\cC))$; the
  value on $S$ is the functor $X \mapsto X^{\otimes S}$. Composing
  with the symmetric monoidal functor $\Free_G \to \Fin^{BG}$ gives a
  symmetric monoidal functor
  $\Free_G \to \Fun_\otimes(\Env(\cC),\Env(\cC))^{BG}$. By
  \cite[III.3.7]{scholze-nikolaus-tc}, this structure is adjoint to a
  symmetric monoidal functor $\Free_G \times_{\Fin} \Env(\cC) \to
  \Env(\cC)^{BG}$.
\end{proof}

\begin{defn}
  For a diagram $(X,f)\co K \to \Env(\cC)_{/BG}$, represented by a map
  $X\co K \to \Env(\cC)$ and a map $f\co |X| \to BG$, we define the
  diagram obtained by \emph{unwinding} $X$ to be the composite
  \[
    \psi^f X\co K \too{(X,f)} \Env(\cC)_{/BG} \too{\psi} \Env(\cC)^{BG}.
  \]
\end{defn}

\begin{prop}
  The composite $\psi^f X\co K \to \Env(\cC)^{BG} \to Fin^{BG}$ is the
  diagram of $G$-sets classified by the map $K \to \Free(G)$. In
  particular, on taking shapes there is a principal $G$-bundle
  $|\psi^f X| \to |X|$, classified by the map $f\co |X| \to BG$.
\end{prop}

\begin{exam}
  There is a canonical prinicipal $C_n$-bundle
  $sd_n S^1 \to S^1 \too{[n]} BC_n$ over the simplicial circle, and
  the unwinding $\psi^{[n]} Z(A)$ of the cyclic bar construction is
  the simplicial subdivision $sd_n Z(A)$
  \cite{bokstedt-hsiang-madsen}. More generally, if $f\co P \to B$ is
  a principal $G$-bundle and $A$ is a commutative algebra then
  $\psi^f (A^{\otimes B}) = A^{\otimes P}$.

  By contrast, the unwinding $\psi^{[2]} Z(A;M)$ of the cyclic bar
  construction with coefficients is a simplicial
  object
  \[
    M \otimes M \Leftarrow M \otimes A \otimes M \otimes A \Lleftarrow
    M \otimes A^{\otimes 2} \otimes M \otimes A^{\otimes 2} \cdots
  \]
\end{exam}

Reorganizing terms, the above can be regarded as the cyclic bar
construction of $A^{\otimes 2}$ with a particular bimodule structure
on $M^{\otimes 2}$, as in the following definition.

\begin{defn}
  \label{def:twistedbimodule}
  Suppose that $A$ is an algebra in a symmetric monoidal
  $\infty$-category $\cC$ and that $M$ is a $k$-linear
  $A$-bimodule. Fix an $n > 0$, and let
  $\tau \co A^{\otimes n} \to A^{\otimes n}$ be a cyclic permutation
  generating an action of $C_n$. The \emph{twisted tensor power
    $M^{\circlearrowleft n}$} is the pullback of the ordinary
  $A^{\otimes n}$ bimodule $M^{\otimes n}$ along the map
  $1 \otimes \tau\co (A^{\otimes n}) \otimes (A^\op)^{\otimes n} \to
  (A^{\otimes n}) \otimes (A^\op)^{\otimes n}$.
\end{defn}

\begin{rmk}
  This twisted bimodule is $C_n$-equivariant with respect to the twist
  maps on $M^{\circlearrowleft n}$ and $A^{\otimes n}$.
\end{rmk}

\begin{prop}
  There is a $C_n$-equivariant natural equivalence of simplicial objects
  \[
    \psi^{[n]} Z(A;M) \simeq Z(A^{\otimes n};
    M^{\circlearrowleft n})
  \]
  in $\Env(\cC)$.
\end{prop}

\begin{rmk}
  These ``cyclic'' versions of $\THH$ with coefficients have also
  appeared in the work of Lindenstrauss--McCarthy
  \cite{lindenstrauss-mccarthy-taylortower} and Malkiewich--Ponto
  \cite{malkiewich-ponto-periodicpoints}.
\end{rmk}

\section{The Tate diagonal}

Fix a cyclic group $C_p$ of prime order. For a based space $W$, there
is a natural space-level diagonal map
\[
  W \to (W^{\wedge p})^{C_p}.
\]
If $X$ is a spectrum, then assembling the space-level diagonal maps
gives a map called the \emph{Tate diagonal}
\[
  X \to (X^{\otimes p})^{tC_p},
\]
constructed by Greenlees--May in \cite{greenlees-may-tate} and
recently developed further in \cite{scholze-nikolaus-tc}. The Tate
diagonal has a number of very useful properties: it is natural in $X$,
it is impervious to the action of $C_p$ on $X^{\otimes p}$, and it is
lax symmetric monoidal. The compatibility between these properties is
expressed as follows.

\begin{thm}[{\cite[III.3.8]{scholze-nikolaus-tc}}]
  For a finite free $C_p$-set $T$ with quotient $\overline T$ and an
  indexed tuple $\{X_{\bar t}\}_{\bar t \in \overline T}$ of spectra,
  there is a Tate diagonal
  \[
    \bigotimes_{\bar t \in \overline T} X_{\bar t} \to \left(\bigotimes_{t \in
        T}X_{\bar t}\right)^{tC_p}.
  \]
  The Tate diagonal is essentially unique as a $BC_p$-equivariant lax
  symmetric monoidal transformation between functors
  $\Free(C_p) \times_{\Fin} \Sp^\otimes_\act \to \Sp$.
\end{thm}

Our notation expresses this in the following way. The Tate diagonal is
a lax symmetric monoidal natural transformation
\[
  \bigotimes^{\Sp} \circ X \to \left(\bigotimes^{\Sp} \circ \psi^f
    X\right)^{tC_p}
\]
defined on $(X,f)$ in $\Env(\Sp)_{/BC_p}$.

\begin{exam}
  The lax symmetric monoidal structure then makes it possible for us
  to study the relationship with module structures. Given a
  commutative ring spectrum $k$, the iterated multiplication map
  $k^{\otimes p} \to k$ is $C_p$-equivariant and so there is a composite
  map
  \[
    \phi\co k \to (k^{\otimes p})^{tC_p} \to k^{tC_p}
  \]
  called the \emph{Tate-valued Frobenius}
  \cite[IV.1.1]{scholze-nikolaus-tc}.\footnote{There are actually two
    maps $k \to k^{tC_p}$ of commutative algebras. One is the canonical
    unit $k \to k^{hC_p} \to k^{tC_p}$ because $k$ has trivial
    $C_p$-action, and the other is $\phi$.}   
\end{exam}

In these terms, we obtain the following indexed Tate diagonal.

\begin{cor}
  \label{cor:nonadjointdiagonal}
  Given a sifted index category $K$, there is a natural lax symmetric
  monoidal natural transformation
  \[
    \hocolim_K \left(\bigotimes^{\Sp} X\right) \to \left(\hocolim_K
      \bigotimes^{\Sp} \psi^f X\right)^{tC_p},
  \]
  of functors $\Fun(K,\Env(\Sp)_{/BC_p}) \to \Sp$.
\end{cor}

\begin{proof}
  When $K$ is sifted, the functor $\hocolim_K$ is lax symmetric
  monoidal by Proposition~\ref{prop:laxhocolim}.
\end{proof}

\begin{exam}
  Suppose $E \to B$ is a principal $C_p$-bundle and that $k$ is a
  commutative ring spectrum. Then the Loday constructions for $B$ and
  $E$ are related by a Tate diagonal:
  \[
    k^{\otimes B} \to \left(k^{\otimes E}\right)^{tC_p}
  \]
\end{exam}

\begin{exam}
  Let $K$ be $\Delta^\op$, the simplicial index category. When applied
  to the cyclic bar construction $Z(A,M)$ in $\Env(\Sp)$, the Tate
  diagonal becomes a natural transformation
  \[
    \THH(A;M) \to \left[\THH(A^{\otimes p}; M^{\circlearrowleft p})\right]^{tC_p}.
  \]
  on Hochschild homology with coefficients.
\end{exam}

\section{A relative Tate diagonal}

The Tate diagonal from Corollary~\ref{cor:nonadjointdiagonal}
takes place in the category of spectra. In this section we will
examine the extent to which this admits a relative version, where $X$
is a diagram of modules over a commutative ring spectrum $k$ and we
attempt to replace the monoidal structure of $\Sp$ with the monoidal
structure in $k$-modules. In our final application, $k$ will be an
Eilenberg--Mac Lane spectrum.

\begin{thm}
  \label{thm:generaldiagonal}
  Suppose that $k$ is a commutative ring spectrum, $K$ is a sifted
  index category, and $(X,f)\co K \to \Env(\LMod_k)_{/BC_p}$ is a
  diagram with shape $|X|$. Then there is a natural \emph{relative Tate
    diagonal}
  \[
    k^{tC_p} \mathop\otimes_{k^{\otimes |X|}} \left(\hocolim_K
      \bigotimes^{\Sp} X\right) \to \left(\hocolim_K
      \bigotimes^{\LMod_k} \psi^f X\right)^{tC_p}.
  \]
  When $k$ is the sphere spectrum, this recovers the ordinary Tate
  diagonal.
\end{thm}

\begin{proof}
  Lax symmetric monoidality implies that the Tate diagonal
  \[
    \varphi\co k^{\otimes |X|} \to \left(k^{\otimes |\psi^f X|} \right)^{tC_p}
  \]
  is a map of commutative ring spectra, and that the Tate diagonal
  \[
    \left(\hocolim_K
      \bigotimes^{\Sp} X\right) \to \left(\hocolim_K
      \bigotimes^{\Sp} \psi^f X\right)^{tC_p}
  \]
  is compatible with the $k^{\otimes |X|}$-module structure on the
  source and the $k^{\otimes |\psi^f X|}$-module structure on the
  target. Similarly, the augmentation map
  \[
    k^{\otimes |\psi^f X|} \to k^{\otimes_k |\psi^f X|} \simeq k
  \]
  is a $C_p$-equivarant map of commutative ring spectra, and the map
  \[
    \hocolim_K \bigotimes^{\Sp} \psi^f X \to \hocolim_K
    \bigotimes^{\LMod_k} \psi^f X
  \]
  is a $C_p$-equivariant map of $k^{\otimes |\psi^f X|}$-modules; we
  can then apply Tate spectra.

  Putting these together, there is a composite map
  \[
    \left(\hocolim_K
      \bigotimes^{\Sp} X\right) \to \left(\hocolim_K
      \bigotimes^{\LMod_k} \psi^f X\right)^{tC_p}.
  \]    
  This is a map of $k^{\otimes |X|}$-modules, with the target module
  pulled back from $k^{tC_p}$. The adjoint map is the desired relative
  Tate diagonal.
\end{proof}

% \todo{also,
%   existence of the map is just functoriality of the Tate diagonal}
\begin{exam}
  Let both $K$ and the map $K \to \Fin$ be trivial. Then the
  $k$-module Tate diagonal
  \[
    k^{tC_p} \otimes_k M \to (M^{\otimes_k p})^{tC_p}
  \]
  is adjoint to the composite
  \[
    M \to (M^{\otimes p})^{tC_p} \to (M^{\otimes_k p})^{tC_p}
  \]
  of the ordinary Tate diagonal with the lax monoidal natural
  transformation on tensor powers. Because it is easy to overlook, it
  is worth noting explicitly that the $k^{tC_p}$-module structure on
  the target makes it into a $k$-module via both the canonical unit
  and the Tate-valued Frobenius; only the Tate-valued Frobenius makes
  this composite a map of $k$-modules.
  
  From the point of view of genuine-equivariant homotopy theory, the
  Tate-valued Frobenius is the structure needed to lift $k$ to a
  $C_p$-equivariant commutative ring spectrum. The $k$-module Tate
  diagonal lifts $M^{\otimes_k p}$ to a $C_p$-equivariant module
  called the relative norm of $M$.
\end{exam}

\begin{exam}
  In the case of relative $\THH$, this becomes a \emph{relative
    $\THH$-diagonal} 
  \[
    k^{tC_p} \otimes_{\THH(k)} \THH(A) \to \left[\THH^k(A)\right]^{tC_p},
  \]
  adjoint to the composite
  \[
    \THH(A) \to  \left[\THH(A)\right]^{tC_p} \to
    \left[\THH^k(A)\right]^{tC_p}
  \]
  of the ordinary $\THH$ diagonal with the base-change map. Again, the
  $k^{tC_p}$-module structure on the target makes it into a
  $\THH(k)$-module via the composite of the augmentation $\THH(k)$
  with the Tate-valued Frobenius $k \to k^{tC_p}$. Similarly, we have
  a relative $\THH$-diagonal with coefficients
  \[
    k^{tC_p} \otimes_{\THH(k)} \THH(A;M) \to \left[\THH^k(A^{\otimes
        p}; M^{\circlearrowleft p})\right]^{tC_p}.
  \]
\end{exam}

% However, if the principal $C_p$-bundle has a chosen trivialization,
% this is automatic. Then there is a $C_p$-equivariant map $|\psi^f X|
% \to C_p$, and naturality of the Tate-valued Frobenius provides a
% commutative diagram
% \[
%   \xymatrix{
%     k^{\otimes |X|} \ar[r] \ar[d]
%     & \left(k^{\otimes |\psi^f X|} \right)^{tC_p} \ar[r] \ar[d]
%     & k^{tC_p} \ar@{=}[d] \\
%     k \ar[r]
%     & \left(k^{\otimes p}\right)^{tC_p} \ar[r] &
%     k^{tC_p}
%   }
% \]
% of commutative ring spectra.

\section{Nonexistence of a true relative diagonal}

We will use a calculation with topological Hochschild homology to
illustrate the nonexistence of a Tate diagonal for $k$-modules,
analogous to the Tate diagonal for spectra. We learned this result
from Lars Hesselholt.

Suppose there was a $k$-module Tate diagonal
\[
  \bigotimes^k_{s \in S} M_s \to \left[\bigotimes^k_{t \in T}
    M_{f(t)}\right]^{tC_p},
\]
compatible with the one for spectra and functorial in pairs of a
principal $C_p$-bundles $f\co T \to S$ and an $S$-indexed tuple of
$k$-modules. Using these compatibilities, we could construct a
diagonal for $\THH$ relative to $k$: a map
\[
  \THH^k(A) \to \THH^k(A)^{tC_p}
\]
for any $k$-algebra $A$, accepting a natural transformation from the
Tate diagonal $\THH(A) \to \THH(A)^{tC_p}$. Specializing to the case
$A = k$, we would get a commutative diagram of ring spectra of the form
\[
  \xymatrix{
    \THH(k)  \ar[r] \ar[d] & \THH(k)^{tC_p} \ar[d] \\
    k \ar@{.>}[r] & k^{tC_p}.
  }
\]

However, Hesselholt--Madsen's calculations in the case of an
Eilenberg--Mac Lane spectrum for $\mb F_p$ (or, more generally, for a
perfectoid ring by work of Bhatt--Morrow--Scholze) show that
this would give a commutative diagram of graded rings
\[
  \xymatrix{
    \mb F_p[u]  \ar[r] \ar[d] & \mb F_p[u^{\pm 1}] \ar[d] \\
    \mb F_p \ar@{.>}[r] & \mb F_p[u^{\pm 1}] \cdot \{1, v\}
  }
\]
upon taking coefficients.

\begin{rmk}
  \label{rmk:krause-nikolaus}
  The construction of this diagram is very close to equivalent to the
  existence of a natural $p$-cyclotomic structure on $\THH^k(A)$
  compatible with that on $\THH(A)$. More explicitly, if $k$ is
  connective and the Tate-valued Frobenius $k \to k^{tC_p}$ can be
  made compatible with the trivial $S^1$-action on the source and the
  natural action on the target, the work of Nikolaus--Scholze lifts
  $k$ to a $p$-cyclotomic ring spectrum \cite{scholze-nikolaus-tc},
  allowing us to take the base-change formula
  \[
    \THH^k(A) \simeq k \otimes_{\THH(k)} \THH(A)
  \]
  and use it to define $\THH^k(A)$ as a $p$-cyclotomic spectrum
  \cite{cyclotomic}. This is possible for the spherical group algebra
  $\mb S[\mb N]$, which is exploited to great effect in recent work of
  Krause--Nikolaus on discrete valuation rings
  \cite{krause-nikolaus-periodicity}.
\end{rmk}

\section{Smooth algebras}
\label{sec:smooth-algebras}

In this section, we assume that $k$ is a commutative ring
spectrum---in particular, $k$ could be the Eilenberg--Mac Lane
spectrum associated to a commutative ring.

\begin{defn}
  Let $A$ be a $k$-algebra and $p$ a prime. We say that a
  $k$-linear $A$-bimodule $M$ \emph{satisfies Tate descent at $p$}
  if the relative $\THH$ diagonal
  \[
    k^{tC_p} \otimes_{\THH(k)} \THH(A,M) \to (\THH^k(A^{\otimes_k p},
    M^{\circlearrowleft_k p})^{tC_p}
  \]
  is an equivalence. If $M$ satisfies Tate descent at all primes, we
  simply say that $M$ satisfies Tate descent.
\end{defn}

\begin{prop}
  \label{prop:thickdescent}
  The collection of $k$-linear $A$-bimodules satisfying Tate descent
  at $p$ is a thick subcategory, and in particular is closed under
  finite limits and colimits.
\end{prop}

\begin{proof}
  The $k$-module $\THH$ diagonal is a natural transformation of exact
  functors: both source and target preserve cofiber sequences. In
  particular, the collection of objects for which the $\THH$ diagonal
  is an equivalence is a thick subcategory of the category of
  $k$-linear $A$-bimodules.
\end{proof}

\begin{prop}
  \label{prop:extendeddescent}
  Any $k$-linear $A$-bimodule of the form $N \otimes_k A$, where $N$
  is a left $A$-module that is perfect as a $k$-module, satisfies Tate
  descent.
\end{prop}

\begin{proof}
  The natural map $k \to A$ induces natural equivalences
  \[
    \THH(k;N) \to \THH(A; N \otimes_k A)
  \]
  and
  \[
    \THH(k^{\otimes p};N^{\circlearrowleft p}) \to \THH(A^{\otimes p};
    (N \otimes_k A)^{\circlearrowleft p}).
  \]
  Therefore, by naturality of the Tate diagonal it suffices to show
  this result when $A = k$. Because the collection of $N$-modules
  satisfying Tate descent is a thick subcategory, it suffices to show
  this result when $N = k$ in order to conclude it is true for all
  perfect $k$-modules.

  In this case, we are considering the relative $\THH$-diagonal
%  \todo{R. says he should have understood what this section is saying
%    concretely, but didn't.}
  \[
    k^{tC_p} \otimes_{\THH(k)} \THH(k) \to \left[\THH^k(k^{\otimes_k p};
      k^{\circlearrowleft_k p})\right]^{tC_p},
  \]
  which simplifies to the natural transformation
  \[
    k^{tC_p} \otimes_{\THH(k)} \THH(k) \to \left[\THH^k(k;
      k)\right]^{tC_p}.
  \]
  Both sides are weakly equivalent to $k^{tC_p}$. Moreover, on both
  sides this equivalence is induced by the map
  \[
    k \to (k^{\otimes p})^{tC_p} \to k^{tC_p}
  \]
  in degree $0$ of the simplicial diagrams defining $\THH$.
\end{proof}

\begin{prop}
  \label{prop:smoothdescent}
  If $A$ is smooth, then all $k$-linear $A$-bimodules which are
  perfect as left $k$-modules satisfy Tate descent.
\end{prop}

This is a reformulation of Theorem~\ref{thm:main}.

\begin{proof}
  Fix any $A$-bimodule $M$ that is perfect over $k$, and let
  $\mathcal{T}$ be the full subcategory of $k$-linear $A$-bimodules
  $B$ such that the bimodule $M \otimes_A B$ satisfies Tate descent. By
  Proposition~\ref{prop:extendeddescent}, the bimodule $A \otimes_k A$
  is in $\mathcal{T}$. The category $\mathcal{T}$ is a thick
  subcategory by Proposition~\ref{prop:thickdescent}. By definition,
  since $A$ is smooth over $k$, $A$ lies in the thick subcategory of
  $k$-linear $A$-bimodules generated by $A \otimes_k A$, and therefore
  $A$ is in $\mathcal{T}$.

  The equivalence of bimodules $M \simeq M \otimes_A A$ then shows
  that $M$ satisfies Tate descent.
\end{proof}

\begin{cor}
  If $A$ is smooth and proper, then $A$ satisfies Tate descent as a
  bimodule over itself.
\end{cor}

\section{Hochschild homology and spectral sequences}

We will now complete proofs of our algebraic statements.

\begin{proof}[Proof of Theorem~\ref{thm:mainhomological}.]
Suppose that $\mb F$ is a perfect field of characteristic $p$, $A$ is a
homologically smooth differential graded $\mb F$-algebra, and
$M$ is a bounded $A$-bimodule. We need to construct a Tate
cohomology spectral sequence
\[
  \tH^*(C_p; HH^{\mb F}_*(A, \underbrace{M \otimes^L_A \dots \otimes^L_A
    M}_{p})) \Rightarrow \tH^*(C_p; \mb F) \otimes_{\mb
    F[u]}\THH_*(A;M).
\]
By the results of \S\ref{sec:hochschild}, the equivalence $\theta$
from $\Ch(\mb F)$ to the category of left modules over the
Eilenberg--Mac Lane spectrum $H\mb F$ takes the totalization of the
Hochschild homology object $\HH^{\mb F}(A,M)$ to the relative topological
Hochschild homology object $\THH^{H\mb F}(\theta A, \theta M)$, and the
cyclic tensor $\HH^{\mb F}(A^{\otimes p},M^{\otimes p})$ to
$\THH^{H\mb F}(\theta A^{\otimes p}, \theta M^{\circlearrowleft
  p})$. Since $A$ is homologically smooth over $\mb F$, $\theta A$ is
smooth over $H\mb F$, and since $M$ is bounded and finitely generated over
$\mb F$, $\theta M$ is a perfect $H\mb
F$-module. Theorem~\ref{thm:main} then implies that the relative
$\THH$ diagonal
\[
  H\mb F^{tC_p} \otimes_{\THH(H\mb F)} \THH(\theta A, \theta M) \to
  \left[\THH^{H\mb F}(\theta A^{\otimes_{H\mb F} p}, \theta
    M^{\circlearrowleft_{H\mb F} p})\right]^{tC_p}
\]
is an equivalence.

The map $\pi_* \THH(H\mb F) \to \pi_* H\mb F^{tC_p}$ is an inclusion
$\mb F[u] \to \mb F[u^{\pm 1}] \cdot \{1, v\}$; in particular, it is a flat
map of rings. Therefore, the K\"unneth formula for the homotopy of a
tensor product identifies the homotopy groups of the source of the
relative $\THH$ diagonal with
\[
  \widehat{H}^*(C_p;\mb F) \otimes_{\mb F[u]} \THH_*(A;M).
\]
On the other hand, the target of the $\THH$ diagonal has a
conditionally convergent Tate spectral sequence
\[
  \widehat{H}^*(C_p; \pi_* \THH^{H\mb F}(\theta A^{\otimes_{H\mb F} p}; \theta
  M^{\circlearrowleft_{H\mb F} p})).
\]
The identification between homotopy and homology under the functor
$\theta$ allows us to re-express this as
\[
  \widehat{H}^*(C_p; \HH_*^{\mb F}(A^{\otimes_{\mb F} p};
  M^{\circlearrowleft_{\mb F} p})).
\]
This Hochschild homology is equivalent to the Hochschild homology of
$A$ with coefficients in $M \otimes^L_A \dots \otimes^L_A M$, as
desired.
\end{proof}

\begin{proof}[Proof of Theorem~\ref{thm:spectralsequence}.]
  Suppose that $\THH_*(A;M)$ is torsion-free as a module over the
  graded ring $\THH_*(\mb F) \cong \mb F[u]$. We wish to show that the
  spectral sequence of Theorem~\ref{thm:mainhomological}, constructed
  in the previous proof, reduces to a spectral sequence with a simpler
  grading.

  The previous theorem gives us a spectral sequence
  \[
    \tH^*(C_p; HH^{\mb F}_*(A, \underbrace{M \otimes^L_A \dots \otimes^L_A
      M}_{p})) \Rightarrow \tH^*(C_p; \mb F) \otimes_{\mb F[u]}\THH_*(A;M).
  \]
  Moreover, both sides are acted on by the Tate cohomology ring
  $\tH^*(C_p;\mb F)$ consisting of permanent cycles, by comparison with $A
  = M = \mb F$.
  
  Since $\THH_*(A;M)$ is graded, torsion-free as a module over
  $\mb F[u]$, and bounded below, we can lift a basis of
  $\THH_*(A;M) / u$ over $\mb F$. An argument by induction on the
  grading shows that this lift is a basis of $\THH_*(A;M)$ over
  $\mb F[u]$: it is isomorphic to a direct sum of shifts of the graded
  module $\mb F[u]$. This freeness of $\THH_*(A,M)$ over $\mb F[u]$,
  together with the base-change equivalence
  \[
    \theta \HH^{\mb F}(A;M) \simeq \THH^{H\mb F}(A;M) \simeq H\mb F
    \otimes_{\THH(\mb F)} \THH(A, M),
  \]
  implies via the K\"unneth spectral sequence that we have an
  isomorphism
  \[
    \HH^{\mb F}_*(A;M) \cong \THH_*(A;M) / (u).
  \]
  Moreover, our chosen lift of basis for $\HH^{\mb F}_*(A;M)$ to
  $\THH_*(A;M)$ gives an isomorphism
  \[
    \THH_*(A;M) \cong \mb F[u] \otimes_{\mb F} \HH_*(A;M)
  \]
  which is noncanonical except up to associated graded.

  If $p=2$, the Tate cohomology ring $\tH^*(C_2;\mb F)$ is a Laurent
  polynomial ring on a generator $x$ satisfying $x^2 = u$. This
  allows us to re-express the Tate spectral sequence in the form
  \[
    \tH^*(C_2; HH^{\mb F}_*(A, M \otimes^L_A M)) \Rightarrow
    \mb F[x^{\pm 1}] \otimes_{k[u]} \THH_*(A;M).
  \]
  As the multiplication-by-$x$ operation makes this spectral sequence
  fully periodic (inducing an isomorphism between degree $d$ and
  degree $(d+1)$ for all $d$), it is equivalent to an \emph{ungraded}
  spectral sequence consisting just of those terms from degree
  zero; this is equivalent to setting $x=1$. This ungraded spectral
  sequence takes the form 
  \[
    \bigoplus_s \tH^s(C_2; HH^{\mb F}_s(A, M \otimes^L_A M)) \Rightarrow
    \mb F \otimes_{\mb F[u]} \THH_*(A;M).
  \]
  Tate cohomology for $C_2$ is isomorphic in all degrees, so we can
  safely replace $\tH^s$ with $\tH^0$. Our noncanonical description
  of $\THH_*(A;M)$ now gives a different identification of the
  right-hand side with $\HH^{\mb F}_*(A;M)$ (now with $u$ sent to $1$,
  rather than $0$).

  At odd primes, the same method applies; however, the Tate cohomology
  ring only contains an inverse of the element $u$ in degree $2$, and
  so setting $u=1$ allows us to interpret it as an even-odd graded
  spectral sequence.
\end{proof}
\bibliography{../masterbib}
\end{document}